\documentclass[11p]{article}

\usepackage{graphicx}%
\usepackage{multirow}%
\usepackage[T1]{fontenc} 
\usepackage{amsmath,amssymb,amsfonts}%
\usepackage{amsthm}%
\usepackage{mathrsfs}%
\usepackage[title]{appendix}%
\usepackage{xcolor}%
\usepackage{textcomp}%
\usepackage{manyfoot}%
\usepackage{booktabs}%
\usepackage{algorithm}%
\usepackage{algorithmicx}%
\usepackage{algpseudocode}%
\usepackage{listings}%
\usepackage{subfig}
\usepackage{pgfplots}
\usepackage{pgfplotstable} % Generates table from .csv
\usepackage{tabularx,longtable,multirow,caption}%hangcaption
\usepackage{booktabs,colortbl}
\pgfplotstableset{
	begin table=\begin{longtable},
		end table=\end{longtable},
}
\usepackage[square,numbers]{natbib}

\newcommand{\e}{\mathbf{e}}
\newcommand{\R}{\mathbb{R}}
\newcommand{\E}{\mathbb{E}}
\newcommand{\sn}{{\mathcal S}^n}
\newcommand{\inprod}[2]{{\langle #1,#2 \rangle}} % inner product
\newcommand{\trace}{\textrm{tr}}
\newcommand{\rank}{\textrm{rank}}
\newcommand{\avg}{\textrm{avg}}
\newcommand{\st}{\textrm{s.t.}}
\newcommand{\diag}{\textrm{diag}}

\usepackage{scalerel}
\usepackage{tikz}
\usetikzlibrary{svg.path}

\definecolor{orcidlogocol}{HTML}{A6CE39}
\tikzset{
	orcidlogo/.pic={
		\fill[orcidlogocol] svg{M256,128c0,70.7-57.3,128-128,128C57.3,256,0,198.7,0,128C0,57.3,57.3,0,128,0C198.7,0,256,57.3,256,128z};
		\fill[white] svg{M86.3,186.2H70.9V79.1h15.4v48.4V186.2z}
		svg{M108.9,79.1h41.6c39.6,0,57,28.3,57,53.6c0,27.5-21.5,53.6-56.8,53.6h-41.8V79.1z M124.3,172.4h24.5c34.9,0,42.9-26.5,42.9-39.7c0-21.5-13.7-39.7-43.7-39.7h-23.7V172.4z}
		svg{M88.7,56.8c0,5.5-4.5,10.1-10.1,10.1c-5.6,0-10.1-4.6-10.1-10.1c0-5.6,4.5-10.1,10.1-10.1C84.2,46.7,88.7,51.3,88.7,56.8z};
	}
}
\newcommand\orcidicon[1]{\href{https://orcid.org/#1}{\mbox{\scalerel*{
				\begin{tikzpicture}[yscale=-1,transform shape]
				\pic{orcidlogo};
				\end{tikzpicture}
			}{|}}}}
\usepackage{hyperref}

\hypersetup{%
  colorlinks=true,
  linkcolor=blue,
  linkbordercolor={0 0 1}
}

\theoremstyle{thmstyleone}%
\theoremstyle{thmstyletwo}%

\theoremstyle{thmstylethree}%

\providecommand{\keywords}[1]
{
  \small	
  \textbf{\textit{Keywords---}} #1
}
\usepackage[multiple]{footmisc}
\makeatletter
\newcommand\footnoteref[1]{\protected@xdef\@thefnmark{\ref{#1}}\@footnotemark}
\makeatother
\title{Ground truth clustering is not the  optimum clustering}
\makeatletter
\newcommand\newtag[2]{#1\def\@currentlabel{#1}\label{#2}}
\makeatother

\date{ }

\begin{document}

\author{Lucia Absalom Bautista \footnote{University of Sevilla, C. San Fernando 4, 41004 Sevilla, Spain \href{mailto:luciixi99@gmail.com}{luciixi99@gmail.com}} \orcidicon{0009-0002-0366-5174}
\and Timotej Hrga \footnote{University of Ljubljana, Faculty of Mathematics and Physics \& Faculty of Mechanical Engineering, Jadranska 19 \& Aškerčeva 6, 1000 Ljubljana, Slovenia~\href{mailto:timotej.hrga@fmf.uni-lj.si}{timotej.hrga@fmf.uni-lj.si}} \orcidicon{0000-0002-4852-1986}
\and Janez  Povh \footnote{University of Ljubljana, Faculty of Mechanical Engineering, Aškerčeva 6, 1000 Ljubljana, Slovenia ;
Rudolfovo – Science and technology center Noov mesto, Podbreznik 15, 8000 Novo Mesto,Slovenia~\href{mailto:janez.povh@rudolfovo.eu}{janez.povh@rudolfovo.eu}} \orcidicon{0000-0002-9856-1476}
\and Shudian Zhao \footnote{Optimization and Systems Theory, Department of Mathematics, KTH Royal Institute of Technology, Lindtstedtsv\"agen 25, 100 44 Stockholm, Sweden~\href{mailto:shudian@kth.se}{shudian@kth.se}} \orcidicon{0000-0001-6352-0968}}
 
% \affil[1]{\orgdiv{}    \orgname{University of Sevilla, C. San Fernando 4}, \city{Sevilla}, \postcode{41004}, \country{Spain}}}

% \affil[2]{ \orgdiv{Faculty of Mathematics and Physics}, \orgname{University of Ljubljana}, \orgaddress{\street{Jadranska 19}, \city{Ljubljana}, \postcode{1000},  \country{Slovenia}}}
% \affil[3]{\orgdiv{Faculty of Mechanical Engineering}, \orgname{University of Ljubljana}, \orgaddress{\street{Aškerčeva 6}, \city{Ljubljana}, \postcode{1000},  \country{Slovenia}}}

% % \affil[3]{\orgdiv{Faculty of Mechanical Engineering}, \orgname{University of Ljubljana}, \orgaddress{\street{Aškerčeva 6}, \city{Ljubljana}, \postcode{1000},  \country{Slovenia}}}
% \affil[4]{ \orgname{Rudolfovo – Science and technology center Noov mesto}, \orgaddress{\street{Podbreznik 15}, \city{Novo Mesto}, \postcode{8000},  \country{Slovenia}}}

% \affil[5]{\orgdiv{Department of Mathematics}, \orgname{KTH - Royal Institute of Technology}, \orgaddress{\street{Lindtstedtsvägen 25}, \city{Stockholm}, \postcode{100 44}, \country{Sweden}}}

%%================================%%
%% Sample for structured abstract %%
%%================================%%

\maketitle
\begin{abstract}
The clustering of data is one of the most important and challenging topics in data science. The minimum sum-of-squares clustering (MSSC) problem asks to cluster the data points into $k$ clusters such that the sum of squared distances between the data points and their cluster centers (centroids) is minimized. This problem is NP-hard, but there exist exact solvers that can solve such problem to optimality for small or medium size instances.  

In this paper, we use a branch-and-bound solver based on semidefinite programming relaxations called SOS-SDP to compute the optimum solutions of the MSSC problem for various $k$ and for multiple datasets, with  real and artificial data, for which the data provider has provided ground truth clustering.

Next, we use several extrinsic and intrinsic measures to evaluate how the optimum clustering and ground truth clustering matches, and how well these clusterings perform with respect to the criteria underlying the intrinsic measures. Our calculations show that the ground truth clusterings are generally far from the optimum solution to the MSSC problem. Moreover, the intrinsic measures evaluated on the ground truth clusterings are generally significantly worse compared to the optimum clusterings. However, when the ground truth clustering is in the form of convex sets, e.g., ellipsoids, that are well separated from each other, the ground truth clustering comes very close to the optimum clustering.

\keywords{minimum sum-of-squares clustering, ground truth clustering, extrinsic measures, intrinsic measures}

\end{abstract}

\section{Introduction}
\label{intro}
\subsection{Motivation}
Data science is strongly related to  the progress in multiple sub-fields of mathematics and computer science. Understanding the data and finding hidden relations between data instances or variables are getting more and more important, especially in the era of Big Data and Artificial Intelligence. Grouping the data instances according to their inner similarity is called  clustering analysis~\cite[Ch.10.3]{GaWiHa:13},~\cite{berkhin2006survey, XuWu:05}. 
A vast amount of literature is devoted to this problem, covering several different aspects. Nevertheless, two important questions  always need to be answered: (i) what is the number of groups that we want to cluster the data into? (ii) what is the similarity measure that will be used to compute the similarity of two data points  and consequently to define groups of similar data, i.e., the clusters?
Usually, these two questions are strongly related, i.e., the number of clusters and the clustering are strongly dependent on the underlying similarity measure.

In some cases, we know in advance the actual number of groups into which the data should be clustered, e.g., based on some domain-specific knowledge. In these cases, the main task is to compute an assignment of the data points into this number of groups, which is also called the $k$-clustering problem.  However, this is a rather rare situation. Usually, the number of clusters has to be determined during the clustering analysis. Some methods, such as the hierarchical clustering methods, determine the number of clusters during the clustering process. Other methods, such as the $k$-means or the $k$-median clustering algorithms, require this number as input, and there are many different criteria and techniques that can be used to determine the appropriate number of clusters, e.g., the cross-validation, the ``Elbow'' method, the gap statistic, the silhouette method, etc., see \cite{hennig2015handbook,kodinariya2013review}.

Clustering  is often part of an exploratory research and in general, there is not only one good/optimum clustering. The selection of similarity measures, number of groups, and good assignment into clusters are often done iteratively, based on an increasing understanding of the data and partial clustering results.
Many algorithms for clustering analysis are known in the literature. Usually, they iteratively define groups of similar data points by following some principles that should hold for good clustering, e.g., cluster homogeneity, cluster completeness, possessing a rag bag cluster, $n$-invariance, see e.g.,~\cite{amigo2009comparison}. In each iteration, they improve the clustering such that the selected criteria  are improved. Once there is no more improvement or the maximum number of iterations is reached, the algorithm terminates. 

In a very complex situation where there are many criteria and principles for good clustering, and consequently many different algorithms for computing such solutions, usually leading to many different clustering solutions for the same dataset, the question naturally arises as to what is a ground truth, i.e., which groups are most natural, which data points really belong to the same group. Very often such clustering is not known, but for some datasets it is provided by the data provider. If we know the ground truth clustering, we can compare the clusterings computed by different algorithms with it. This was the main motivation for us. We focused on  an exact mathematical programming formulation of the clustering problem, solved it with available global optimization algorithms, and compared the optimum solution with the ground truth clustering.

%Many data scientists implement ad hoc clustering algorithms and then evaluate them by generic measures.??
\subsection{The data clustering problem}

The data clustering problem can be formulated as follows: given a  set of data points $\mathcal{P}=\{p_1,p_2,\ldots,p_n\}\subset \mathbb{R}^m$,  the objective is to find the cluster number $k$ and an assignment $\chi \colon \mathcal{P}\rightarrow \{1,2,\ldots,k\}$  such that $\chi(p_i)=\chi(p_j)$ for any pair of (very) data similar points $p_i$ and $p_j$ ($p_i \neq p_j$). For any $k$ assignment $\chi$ we can define the $i$-th cluster, implied by $\chi$ as $C_i=\chi^{-1}(i)$ and the corresponding clustering is 
$\mathcal{C}=\{C_1,\dots, C_k \}$.

Even if we already know the number $k$, there are still many possibilities for how to measure (high) similarity between the data points and how to compute good clustering $\mathcal{C}$.
In this paper, we will focus on  Euclidean distance measure between the two data points. The closer the points are in the Euclidean distance, the more similar they are. We  evaluate each clustering  by considering the total sum of squares of distances  between the data points and the centroid of the corresponding cluster. This means that : (i) for each cluster we compute its centroid (the vector sum of all data points from this cluster, divided by the cardinality of this cluster). (ii)  for each data point we compute the square Euclidean distance between it and the centroid of the cluster to which this data point is assigned.
The sum of squares of these distances across all data points is a measure of the quality of the clustering. The smaller it is, the better the  clustering.

We can   formulate the problem of minimizing the sum of squares of distances as a mathematical optimization problem, where the decision variables represent  assignments $\chi$ and the objective function is 
the sum of squares of distances. This yields the so-called 
 minimum sum-of-squares clustering (MSSC) problem, formally introduced in Section \ref{subsec:SOS-SDP}.
This problem is  well-known in the literature \cite{rao1971cluster,vinod1969integer}
as it arises in a wide range of applications, for example, image segmentation \cite{dhanachandra2015image, shi2000normalized}, biology \cite{jiang2004cluster} and document clustering \cite{mahdavi2009harmony}. 
It is  an NP-hard problem \cite{MAHAJAN201213,megiddo1984complexity}, which means that there is no polynomial time algorithm to solve it to optimality  (unless P=NP). Nevertheless,  we can still solve it for small datasets, e.g., by using the exact solver SOS-SDP \cite{piccialli2022sos}.
The announcement of this solver was actually the trigger  for our research. We decided to solve the MSSC problem to optimality for a number of small datasets for which the ground truth clustering $\mathcal{C}_{true}$ and the true number of clusters $k_{true}$ were known - usually provided by the dataset provider.

Clustering performance can be quantified by a number of metrics.
There are mainly two types of metrics to evaluate this performance.
To compare different clustering methods, we use intrinsic and extrinsic measures. Extrinsic measures require ground truth labels, while intrinsic measures do not.

The datasets with a ground truth clustering are actually the  classification datasets, where the classification labels are considered as  the ground truth label. This idea was elaborated in  \cite{zhang2019critical}, where the authors show that this approach requires careful attention since the class labels are assigned based on the properties of each individual data point, while the clusterings take into account relations between the data. Our paper is aligned with this observation and provides a deeper understanding of this issue.

A ground truth-based comparative study was done in \cite{zhu2008ground}, where the authors did not solve an exact clustering problem (like MSSC) but rather compared five widely known approximate algorithms for data clustering on  seven published microarray gene expression datasets and one synthetic 
dataset. The  performances of these algorithms were assessed with several  quantitative performance measures, which are different than the measures that we use in this paper.

In \cite{gribel2019hg} the authors introduced a population-based metaheuristic algorithm to solve MSSC approximately, which  performs like  a multi-start $k$-means.
They demonstrate that this algorithm outperforms  all recent state-of-the-art approximate algorithms for MSSC in terms of  local minima and the computed  clusters  are  closer to the ground truth compared to the clusters computed by other algorithms, on the synthetic Gaussian-mixture datasets.

In  \cite{amigo2009comparison,zhang2019critical,zhu2008ground}, the ground truth labels and the extrinsic metrics were used to validate the quality of clusterings and compare different clustering methods. 
The ground truth labels have  shown limitations in many instances for detecting the hidden pattern of the sample and they can produce misleading information when used as true ground labels for measuring the quality of clusterings \cite{zhang2019critical}.

\subsection{Our contributions}
The main goal of this paper is to solve exactly the MSSC problem on a large list of small or medium size datasets available in the literature for which the true number of clusters  $k_{true}$  and the  ground truth clustering $\mathcal{C}_{true}$ are known, and to check how the MSSC optimum clustering $\mathcal{C}_{MSSC}$, i.e. the clustering that corresponds to the optimum solution of MSSC aligns with the ground truth clustering. 
More precisely, we
\begin{itemize}
    \item solve the mathematical programming formulation of clustering problem MSSC to optimality for  values of $k$  with
    $|k-k_{true}|\le 2$ for 12 real datasets and for 12 artificial datasets, obtained from \url{https://github.com/deric/clustering-benchmark/tree/master/src/main/resources/datasets}, which are small \\ enough and for which the ground truth clustering is available;
    \item  compare  the MSSC optimum clusterings $ \mathcal{C}_{opt}$ with the ground truth clusterings $\mathcal{C}_{true}$ by computing a number of extrinsic measures: Adjusted Mutual Information (AMI), Adjusted Random Score (ARS), Homogeneity (h), Completeness (c), normalized mutual information (NMI) and  Fowlkes-Mallows scores (FMS);
    \item  additionally evaluate the quality of the optimum MSSC clusterings $\mathcal{C}_{MSSC}$ and of the ground truth clusterings  $\mathcal{C}_{true}$ by computing three intrinsic measures: the Calinski-Harabasz Criterion  (CHC), the Davies Bouldin Index  (DBI), and    the Silhouette Evaluation Score $S_{score}$; 
    \item provide a better understanding of why and when  the optimum clustering aligns with the ground truth clustering. 
\end{itemize}
Our analysis shows that
\begin{itemize}
    \item the optimum clustering (i.e., the optimum solution of MSSC) can be computed if the number of data points times the number of clusters $k$  is approximately below  $1000$;
    \item the  optimum  clusterings  $\mathcal{C}_{MSSC}$ usually  significantly differ from the ground truth clustering in the following manner: the values of the  extrinsic measures  are often far below the optimum, which is 1 and would be achieved if the ground truth and the optimum clustering would be the same. Additionally, the values of extrinsic measures, evaluated at $\mathcal{C}_{MSSC}$ and corresponding to $k_{true}$ are rarely  optimum, i.e.,  other $k$ often gives better values of these measures;
    \item 
    the ground truth clusterings are  usually much worse compared to the  optimum clusterings, by considering the three  intrinsic quality measures;  
    \item  when the ground truth clustering has natural expected  geometry, i.e.,  the clusters have the form of convex sets, e.g., ellipsoids, which are well separated from each other, then the ground truth clustering is  very similar to the optimum clustering.
\end{itemize}
 
These observations are not unexpected, but to the best of our knowledge have never been shown so explicitly.

\subsubsection{Notation}
 The Euclidean distance between two points $p,q\in\mathbb{R}^m$ is denoted by $d_E$ and defined as follows  $d_E(p,q)=\|p-q\|=\sqrt{\sum_{i}^m (p_i-q_i)^2}$. We use $[n]$ to denote the set of integers $\{1,\dots,n\}$. The trace of matrix $X$ is denoted by $\trace(X)$. 
The space of symmetric matrices is equipped with the trace inner product, which for any $X, Y \in {\mathcal S}^{n}$ is defined as $\inprod{X}{Y}:= \trace (XY)$. The associated norm is the Frobenius norm  $\| X\|_F := \sqrt{\trace (X^2)}=\sqrt{\sum_{ij} x_{ij}^2}$.
The cone of symmetric positive semidefinite matrices  of order $n$ is denoted by ${\mathcal S}_+^n :=\{X \in  {\mathcal S}^n\mid X\succeq \mathbf{0} \}$.

We denote by $\e_n$ the vector of all ones of length $n$. 
 In case that the dimension of ${\mathbf e}_n$
is clear from the context, we omit the subscript.
The operator $\diag \colon \R^{n\times n} \to \R^n$ maps a square  matrix to a vector consisting of its diagonal elements. Its adjoint operator is denoted by $\textrm{Diag}\colon \R^n \to \R^{n\times n}$.  The rank of matrix $X$ is denoted by $\rank(X)$. The average value of components of vector $v$ is denoted by $\avg(v)$.

\section{Mathematical programming formulation for MSSC problem}
 \label{subsec:SOS-SDP}

In this paper, we use  the formulation of the MSSC  problem as a mathematical optimization problem in binary variables, which are subject to linear constraints, and with a non-convex objective function that represents the sum of squares of Euclidean distances between the data points and the centroids of the clusters to which the points correspond. 
 
 For given set of data points $\mathcal{P}=\{p_1,p_2,\ldots,p_n\}\subset \mathbb{R}^m$ and given integer $k$ we define  $k$-clustering as  an assignment  $\chi \colon \mathcal{P}\rightarrow \{1,2,\ldots,k\}$.
 The data points that are mapped to integer $1\le i\le k$ are called the $i$-th cluster.
 The assignment is nontrivial if there is no empty cluster.
 We can represent each nontrivial assignment  $\chi$ by a matrix $X\in \{0,1\}^{n\times k}$ such that $x_{ij}=1$ if and only if 
 $\chi(p_i)=j$. Therefore, the row-sums of $X$ must be equal to 1 (each data point is assigned to exactly one cluster)  and the column-sums of $X$ must be at least one (otherwise the assignment is not nontrivial).
 Therefore, the minimum sum-of-squares clustering (MSSC) problem for a fixed $k$ can be formulated as a non-linear integer programming problem (Vinod \cite{vinod1969integer}, Rao \cite{rao1971cluster}):
\begin{equation}\label{eq:mssc}\tag{$\mathrm{MSSC}$}
    \begin{aligned}
    \min~&\sum^{n}_{i=1}\sum^k_{j=1}x_{ij}\|p_i - c_{j}\|^2\\
\st~~& \sum^k_{j=1}x_{ij} = 1, \forall i \in [n],\\
& \sum^n_{i=1}x_{ij} \geq 1, \forall j \in [k],\\
& x_{ij} \in \{0, 1\}, \forall i \in [n] ~\forall j \in [k], \\
& c_j \in \R^m, \forall j \in [k],
    \end{aligned}
\end{equation}
where $c^j$ is the centroid of the $j$-th cluster and can be substituted as $c^j = \frac{\sum^n_{i=1} x_{ij} p_i}{\sum^n_{i=1} x_{ij}}, \forall j \in [k]$. 
% Hence, we can reformulate \eqref{eq:mssc} without $c_j$. 
Peng et al.~\cite{peng2007approximating} introduced equivalent formulations for \eqref{eq:mssc} by introducing substitution $Z:=X(X^\top X)^{-1}X^\top$:
\begin{equation}\label{eq:fullrank-mssc}
    \begin{aligned}
        \min~& -\langle W, Z \rangle\\
        \st~~& Ze = e, \\
            & \trace(Z) = k,\\
             &Z \succeq 0,~Z\ge 0,\\
             &\mathrm{rank}(Z)=k,
    \end{aligned}
\end{equation}
where $W = (W_{ij}) \in \sn$ with $W_{ij} = p_i^\top p_j$.
If we know an optimum matrix $Z$ for \eqref{eq:fullrank-mssc}, we can get back the optimum clustering matrix $X$  for \eqref{eq:mssc} by using the fact that $z_{ij}>0$ if and only if the data vertices $p_i$ and $p_j$ are in the same cluster. So we can put $p_1$ to cluster 1 and the same we do with all $p_i$ with $z_{1i}>0$. The first data point that remains unassigned is then put to cluster 2 and the same all data points with $z_{2i}>0$.
This is repeated until all data points are assigned.

MSSC in formulation \eqref{eq:fullrank-mssc}  remains  NP-hard. However, this formulation opens new possibilities to solve \eqref{eq:mssc}  to optimality. Piccialli et al~\cite{piccialli2022sos} introduced an SDP relaxation for \eqref{eq:fullrank-mssc} by eliminating the constraint $\mathrm{rank}(Z)=k$,  which can be solved to arbitrary precision with state-of-the-art  SDP solvers (e.g., \cite{mosek,sun2020sdpnal+}). SDP relaxation serves as a model to generate lower bounds in the exact methods such as the branch-and-bound algorithm. 
 The SOS-SDP algorithm from \cite{piccialli2022sos} is an implementation  of such   branch-and-bound algorithm based on this relaxation, which was additionally significantly strengthened with the following cutting planes: 

\begin{itemize}
    \item \textbf{Triangle inequalities}, which are based on the observation that if the
points $i$ and $j$ are in the same cluster and the points $j$ and $h$ are in the same cluster, then the points $i$ and $h$ must necessarily belong to the same cluster:
$$
Z_{ij} + Z_{ih} \leq Z_{ii} + Z_{jh}, ~\forall i, j, h \in [n],~ i\neq j\neq h.
$$
\item \textbf{Pair inequalities} 
$$
Z_{ij} \leq Z_{ii},~Z_{ij} \leq Z_{jj},~\forall i, j \in [n],~i \neq j
$$
that every feasible solution of \eqref{eq:fullrank-mssc} satisfies. 
\item \textbf{Clique inequalities}
$$
\sum_{(i,j)\in I\times I,i<j} Z_{ij} \geq \frac{1}{n-k+1},~\forall I \subset [n],~|I| =k+1,
$$
enforcing that for any subset $I$ of $k + 1$ points at
least two points have to be in the same cluster.%, and the number of clusters is $k$. 
\end{itemize}

The branch-and-bound algorithm also demands strong upper bounds for the optimum value, which are usually obtained by computing good feasible solutions, in our case good assignments to clusters. The SOS-SDP algorithm computes these bounds by using the COP $k$-means algorithm from \cite{wagstaff2001constrained}, which is a special variant of the well-known $k$-means heuristics.

\section{Quality measures}\label{sec:quality_measures}
Let denote by $\mathcal{C}_{true}:=\{C_{true,1},\dots,C_{true,k} \}$ the clustering with cluster number $k$, given by the data provider  - we call it the ground truth clustering, and  $\mathcal{C}_{MSSC}=\{C_{MSSC,1},\dots,C_{MSSC,k} \}$  the exact clustering, i.e., the optimum solution of \eqref{eq:mssc},  computed in our case by the SOS-SDP algorithm.

To assess the quality of the computed clusterings, we use both external (extrinsic) and internal (intrinsic) measures. Using the former, we will compare the computed clusters $\mathcal{C}_{MSSC}$ with the ground truth clusters
$\mathcal{C}_{true}$, while the latter do not require ground truth clustering because  for the computed clusters they measure certain criteria such as cluster compactness. There is a long list of possible measures, see for example \cite{amigo2009comparison}, but we restrict ourselves only to those implemented in the Python library {\tt scikit-learn}\footnote{https://scikit-learn.org/stable/modules/clustering.html\#clustering-performance-evaluation}, which we chose for our computations, see also \cite{pedregosa2011scikit}.

\subsection{Extrinsic measures} \label{subsec:methods_for_numer_of_clusters}
The extrinsic methods in general evaluate how two different clusterings match. %Let denote by  $\mathcal{C}:=\{C_1,\dots, C_k \}$ the clustering with cluster number $k$. 
In our case, we use them to assess how the clustering computed by solving the MSSC to optimum (we call this clustering the optimum clustering and denote  it by $\mathcal{C}_{MSSC}$)  matches with the ground truth clustering, i.e., the clustering provided by the data provider algorithm (we denote it by  $\mathcal{C}_{true}$). 
We use six measures from the Python scikit-learn package  \cite{scikit-2011}. For the sake of completeness, we provide here a short description for each of them,  based on the   well-known literature \cite{davies1979cluster,fano1949transmission,fowlkes1983method,hubert1985comparing,rosenberg2007v,ROUSSEEUW198753,shannon2001mathematical}.

\begin{itemize}

      \item {\bf Mutual Information}:
    The mutual information \cite{fano1949transmission,shannon2001mathematical} is also known as the information gain and is computed by
    \begin{equation}
    \begin{aligned}
        \mathrm{ MI}(&\mathcal{C}_{true},\mathcal{C}_{MSSC})=\\
        &\sum_{C_{true,i} \in \mathcal{C}_{true}} \sum_{C_{MSSC,j} \in \mathcal{C}_{MSSC}} \frac{|C_{true,i}\cap C_{MSSC,j}|}{(\sum_{C_{true,i} \in \mathcal{C}_{true}} |C_{true,i}|)^2 } \log\frac{|C_{true,i} \cap C_{MSSC,j}|}{|C_{true,i}||C_{MSSC,j}|},   
    \end{aligned} 
    \end{equation} 
    where $k_{true}$ and $k_{MSSC}$ are the number of clusters in  $\mathcal{C}_{true}$ and $\mathcal{C}_{MSSC}$, respectively, and $C_{\cdot,i}$ is the $i$-th cluster in the clustering assignment $\mathcal{C}_{\cdot}$. The score is nonnegative, and a higher value indicates a higher similarity between the two clusterings. 
    \item {\bf Adjusted Mutual Information (AMI)}:
   
    % \textcolor{red}{Why AMI and NMI in R not equal to AMI and NMI in Python}
    % \todoshudian{check EMI for random clusters}
    The adjusted mutual information score for  $\mathcal{C}_{true}$ and  $\mathcal{C}_{MSSC}$ is computed as~\cite{ami2009}
    \begin{equation*}
    \mathrm{AMI}(\mathcal{C}_{true}, \mathcal{C}_{MSSC}) = \frac{\mathrm{MI}(\mathcal{C}_{true}, \mathcal{C}_{MSSC}) - \E[\mathrm{ MI}]}{
    \avg(H(\mathcal{C}_{true}), H(\mathcal{C}_{MSSC})) - \E[\mathrm{ MI}]},
    \end{equation*}
    where $\avg(\cdot)$ denotes the arithmetic average and $H(\mathcal{C})$ is the entropy of a clustering $\mathcal{C}$:
    \begin{equation}
       H(\mathcal{C}) = \sum_{C_i \in \mathcal{C}} -\frac{|C_i|}{\sum_{C_j \in \mathcal{C}} |C_j|} \log \frac{|C_i|}{\sum_{C_j \in \mathcal{C}} |C_j|},
    \end{equation}
    and $\mathbb{E}(\mathrm{ MI})$ is the expected mutual information between two random clusterings.
    
   The AMI score is between $-1$ and 1. It is 1 when the compared clusterings are identical.

    \item {\bf Adjusted Random Score (ARS)}: ARS \cite{hubert1985comparing,yeung2001details}  computes  the proportion  of pairs of data points that are in the same cluster or in a different cluster in both clusterings (in our case: $\mathcal{C}_{true}$ and $\mathcal{C}_{MSSC}$) and normalize this proportion by using the  expected similarity  specified by a random model, usually based on the generalized hypergeometric distribution. It is computed as
    \begin{equation*}
         \mathrm{ARS} = \frac{\mathrm{ RI} - \E[\mathrm{ RI}]}{\max(\mathrm{ RI}) - \E[\mathrm{ RI}]},
    \end{equation*}
    where random index  $ \mathrm{ RI}= \frac{a+b}{\binom{n}{2}}$ and  $a$ and $b$ are the numbers of pairs of elements that are in the same cluster in $\mathcal{C}_{true}$ and in the same cluster in $\mathcal{C}_{MSSC}$ and  the number of pairs of elements that are in different clusters in $\mathcal{C}_{true}$ and in different clusters in $\mathcal{C}_{MSSC}$, respectively. The value $\E[\mathrm{ RI}]$ is the expected value of $\mathrm{ RI}$ between  two random clusterings.
    The score is between $-1$ and 1. It is 1 when the compared clusterings are identical.

    \item {\bf Homogeneity} ($h$), {\bf Completeness} ($c$) and {\bf V-measure} ($v$)  \cite{rosenberg2007v}:
        % \item Homogeneity \cite{rosenberg2007v}: 
        % A perfectly homogeneous clustering is one where each cluster has data-points belonging to the same class label.
        
        Homogeneity measures how homogeneous the clusters are in $\mathcal{C}_{MSSC}$, i.e., whether the data points from the same cluster from  $\mathcal{C}_{MSSC}$  (mostly) belong to the same cluster in $\mathcal{C}_{true}$.\\
        Completeness measures how many similar samples are put together by the clustering algorithm, i.e., if the data points belonging to the same cluster from $\mathcal{C}_{true}$ are also in the same cluster in $\mathcal{C}_{MSSC}$.
        
        The homogeneity $h$ and completeness $c$ are between 0 and 1 and are defined as
         \begin{equation*}
    h= 1- \frac{H(\mathcal{C}_{true}|\mathcal{C}_{MSSC})}{H(\mathcal{C}_{true})},
    ~c= 1- \frac{H(\mathcal{C}_{true}|\mathcal{C}_{MSSC})}{H(\mathcal{C}_{MSSC})},
         \end{equation*}
         where $H(\mathcal{C}_1\mid \mathcal{C}_{2})$ is the conditional entropy of the clusters in $\mathcal{C}_1$  given the clustering prediction $\mathcal{C}_2$ (recall, $n$ is the number of data points)
         \begin{equation}
             H(\mathcal{C}_{1} \mid \mathcal{C}_{2})= - \sum_{C_{2,j}\in \mathcal{C}_{2}} \sum_{C_{1,i} \in \mathcal{C}_{1}} \frac{|C_{2,j}\cap C_{1,i}|}{n} \log \frac{|C_{2,j}\cap C_{1,i}|}{\sum_{C_{1,s} \in \mathcal{C}_{1}} |C_{2,j}\cap C_{1,s}|},
         \end{equation}
         and $H$
 is entropy, defined above.

        The V-measure $v$ is the weighted harmonic mean between the homogeneity and the completeness: 
    \begin{equation}
       v = \frac{(1 + \beta) \cdot h  \cdot c}{\beta \cdot h + c}, 
    \end{equation}
         where $\beta \geq 0$. When $\beta<1$, more weights are attributed to homogeneity; when $\beta >1$ more weights are attributed to completeness; when $\beta=1$, the score is also called {\bf normalized mutual information (NMI)}, which we report in tables of Section~\ref{sec:results}:
         \begin{equation}
       \mathrm{NMI} = \frac{2 \cdot h  \cdot c}{ h + c}, 
    \end{equation}

    \item {\bf  Fowlkes-Mallows scores (FMS)} \cite{fowlkes1983method}:
    The Fowlkes-Mallows score (FMS) is defined as the geometric mean between pairwise  precision and recall, using True Positive (TP), False Positive (FP), and False Negative (FN):
\begin{equation*}
   \mathrm{FMS} = \frac{\mathrm{TP}}{ \sqrt{(\mathrm{TP + FP}) \cdot (\mathrm{TP + FN})}},    
\end{equation*}
where TP is the number of pairs of points that belong to the same clusters in both $\mathcal{C}_{true}$ and $\mathcal{C}_{MSSC}$, FP is  the number of pairs of points that belong to  the same clusters in $\mathcal{C}_{true}$ and to different clusters in  $\mathcal{C}_{MSSC}$, and FN is the number of  pairs of points that belongs to  the same clusters in $\mathcal{C}_{MSSC}$ and to different clusters in $\mathcal{C}_{true}$. The score ranges from 0 to 1, and a high value indicates a good similarity between the two clusters.

\end{itemize}

\subsection{Intrinsic measures}

Available methods in Python:
\begin{itemize} 
 
    \item Calinski-Harabasz Criterion~\cite{harabasz1974dendrite} (CHC), also known as the Variance Ratio Criterion.

The score is defined as the ratio of the  between-the-cluster dispersion and the within-the-cluster dispersion. Good clustering has a large variance between the cluster  and a small variance within the cluster, hence a large CHC.
    \item Davies Bouldin Index~\cite{davies1979cluster} (DBI), which is defined as follows.
    Let $C_{r}, r\in [k]$, be the $r$-th cluster and  $s_r$  the average Euclidean distance between the data points $p_i$ from cluster $C_r$   and the centroid of that cluster, denoted by $c^r$:
    \begin{equation}
        s_r = (\frac{1}{|C_r|}\sum_{i\in C_r} \|p_i -c^r\|_2^2)^{1/2},\forall r\in [k].
    \end{equation}
The DBI score is defined as
    \begin{equation}
       \mathrm{DBI} = \frac{1}{k}\sum_{i=1}^k\max_{i\neq j \in [k]}(\frac{s_i+s_j}{\|c^i-c^j\|_2}), 
    \end{equation}
Therefore, DBI can only have non-negative values, and the smaller the value of DBI, the better the clustering.  
    
       \item Silhouette Evaluation Score    ~\cite{ROUSSEEUW198753} ($S_{score}$), defined as 
    \begin{equation}
        S_{score} = \frac{1}{n}\sum_{i=1}^n\frac{b_i-a_i}{\max \{a_i,b_i\}},
    \end{equation}
 where $a_i$ is the mean distance between the $i$-th data point and the other data points of the same cluster.
%  associated with some cluster $C_r$ and the other data points of the cluster $C_r$.
 Similarly, $b_i$ is the mean distance between the $i$-th data point and all other data points of the nearest cluster. The fraction 
 $\frac{b_i-a_i}{\max \{a_i,b_i\}}$ thus measures how well the $i$-th data point is positioned in the cluster to which it is assigned. Thus, the value of $S_{score}$ is between $-1$ and $1$, and the closer it is to $1$, the  more appropriate the clustering.
    
\end{itemize}
\section{Numerical Results}\label{sec:results}
\subsection{Datasets}
As benchmark datasets, we used datasets collected by Tomas Barton and available in GitHub \footnote{https://github.com/deric/clustering-benchmark/tree/master/src/main/resources/datasets}.
These datasets consist of a subset of real data and a subset of artificial data. For all these datasets, the data provider has also provided the so-called ground truth clustering, i.e., the number of clusters and the classification into this number of clusters according to certain similarity criteria that include the true groups. We were not aware of this similarity criterion, probably it cannot be  described explicitly, so we used the ground truth labeling only as an input for our research.

From both datasets, we selected only those for which we were capable of solving \eqref{eq:mssc} for the values of $k\ge 2$ with $|k-k_{true}|\le 2$.
Therefore, we exclude the instances where the number of data points or  $k_{true}$ was too large (i.e., $n\cdot k_{true}>1000$).
Nevertheless, we tried to solve also few instances with  $n\cdot k_{true}>1000$ and the computing times were huge (more than a week on one computing node with 2x AMD EPYC 7402 24-core processors and 
128 GB DDR4-3200 ram), but we managed to solve 
five such instances: {\tt dermatology, ecoli, glass, 3MC,  lsun. }

The information about the selected datasets is summarized in Table~\ref{tab:real_data_sum} and Table~\ref{tab:artificial_data_sum}, where $n$ is the number of data points, $m$ is the number of variables (features), including the categorical variable containing the labels of the ground truth clustering in the data samples,  and $k_{true}$ is the clustering number for the ground truth clustering.

All these datasets are formally numerical, which means that there is no categorical value in any dataset. However, a closer look reveals that some real datasets include also  variables that are by their nature categorical, like  {\tt heart-statlog} and {\tt zoo}, which contain half or even all  variables that are categorical  by their nature, respectively. We did not take this fact into account.  The exact clustering results are included in the GitHub repository, see \href{https://github.com/shudianzhao/GroundTrueCluster_VS\_ExactCluster}{https://github.com/shudianzhao/GroundTrueCluster\_VS\_ExactCluster}.

\pgfplotstableread[col sep=comma]{real_data_summary.tex}\data
	\pgfplotstabletypeset[
	font=\footnotesize,
	multicolumn names, % allows to have multicolumn names
	columns={name,type,n,m,K,M},
	fixed zerofill,
	fixed,
 	columns/name/.style={string type,
		column name={instances}
	},
	columns/type/.style={string type,
		column name={type}
	},
	columns/m/.style={int detect,
		column name={$m$},
	},
    columns/n/.style={int detect,
		column name={$n$},
	},
	columns/K/.style={int detect,
		column name={$k_{true}$},
	},
	columns/M/.style={string type,
		column name={measuring type},
	}, 
	    column type=r,
	every head row/.style={
	before row={
	\caption{Summarized information for real data} \label{tab:real_data_sum} \\
	\toprule
% 	&  &  \multicolumn{6}{c|}{extrinsic measurements}& \multicolumn{3}{c}{intrinsic measurements}\\
	},
	after row={
% 	&  &  &  & (rounded) &  (\si{\second}) & (rounded) & (\si{\second}) & \% &  &  &  
	\endfirsthead
	\midrule}
	},
	every last row/.style={after row=\bottomrule}
	]{\data}

\pgfplotstableread[col sep=comma]{artificial_data_summary.tex}\data
	\pgfplotstabletypeset[
	font=\footnotesize,
	multicolumn names, % allows to have multicolumn names
	columns={name,type,n,m,K,M},
	fixed zerofill,
	fixed,
 	columns/name/.style={string type,
		column name={instances}
	},
	columns/type/.style={string type,
		column name={type}
	},
	columns/m/.style={int detect,
		column name={$m$},
	},
    columns/n/.style={int detect,
		column name={$n$},
	},
	columns/K/.style={int detect,
		column name={$k_{true}$},
	},
        columns/M/.style={string type,
		column name={measuring type},
	},  
	column type=r,
	every head row/.style={
	before row={
	\caption{Summarized information for artificial data} \label{tab:artificial_data_sum} \\
	\toprule
% 	&  &  \multicolumn{6}{c|}{extrinsic measurements}& \multicolumn{3}{c}{intrinsic measurements}\\
	},
	after row={
% 	&  &  &  & (rounded) &  (\si{\second}) & (rounded) & (\si{\second}) & \% &  &  &  
	\endfirsthead
	\midrule}
	},
	every last row/.style={after row=\bottomrule}
	]{\data} % filename/path to file

\subsection{Results on real and synthetic datasets}

Table~\ref{tab:real_data_eval_clust} and Table~\ref{tab:artificial_data_eval_clust} contain the names of datasets (column instances) used in the computations, the value $k$ for which \eqref{eq:mssc} was solved,  the corresponding  optimum values $d_{SOS\cdot}$ for \eqref{eq:mssc} (note that the rows highlighted  in gray contain the feasible values of \eqref{eq:mssc}, evaluated at $\mathcal{C}_{true}$). The rest  of the columns contain the scores of the extrinsic and the intrinsic measures, described in Section \ref{sec:quality_measures}, which were computed for the optimum clusterings $\mathcal{C}_{opt}$ and for $\mathcal{C}_{opt}$ (gray rows).

Column $d_{SOS\cdot}$ show that  optimum values of \eqref{eq:mssc} decrease with $k$, as is expected - more groups enable smaller sum of squares of distances between the data points and the cluster centroids. For each dataset, we can observe that the corresponding value of $d_{SOS\cdot}$ in the gray rows, i.e., the values of \eqref{eq:mssc} evaluated at $\mathcal{C}_{opt}$, is  usually much higher compared to the optimum value for  $k=k_{true}$. This reveals that the ground truth clusterings $\mathcal{C}_{true}$ are often far from the optimum of \eqref{eq:mssc}. 

The best scores for each measure, evaluated on optimum clusterings, are highlighted in bold. All the measures  we use  have different preferences among these clustering results: different measures often achieve their best value for different $k$. 

Moreover, the extrinsic measures rarely achieve the best score for $k=k_{true}$ (recall, $k_{true}$ is $k$ in the gray rows) and   the ground truth clusterings  have scores for the intrinsic measures  which significantly differ from the best  scores, for both real  and artificial data set.

% \resizebox{\linewidth}{!}{	
\pgfplotstableread[col sep=comma]{real_data_summary_all_label_clust_eval.tex}\data
 \pgfplotstabletypeset[
	font=\footnotesize,
	multicolumn names, % allows to have multicolumn names
	columns={dataname,k,opt_sos_sdp,adjusted_mutual_info_score,adjusted_rand_score,homogeneity_score,completeness_score,normalized_mutual_info_score,fowlkes_mallows_score,calinski_harabasz_score,davies_bouldin_score,silhouette_score},
	fixed zerofill,
	fixed,
 	columns/dataname/.style={string type,
		column name={instances}
	},
	columns/k/.style={int detect,
		column name={$k$},
		% assign column name/.style={
		% /pgfplots/table/column name={\multicolumn{1}{c|}{##1}}}
	},  
    columns/opt_sos_sdp/.style={precision=0,
        column name ={$d_{SOS}$.},
        assign column name/.style={/pgfplots/table/column name={\multicolumn{1}{c|}{##1}}}
    },
	columns/completeness_score/.style={precision=2,column name={$c$},
	},
    columns/homogeneity_score/.style={precision=2,
	column name={$h$},	},
	columns/normalized_mutual_info_score/.style={precision=2,
	column name={NMI},
	},
	columns/fowlkes_mallows_score/.style={precision=2,
	column name={FMS},
	  assign column name/.style={
		/pgfplots/table/column name={\multicolumn{1}{c|}{##1}}}
	},
	columns/adjusted_mutual_info_score/.style={precision=2,column name={AMI},
	},
	columns/adjusted_rand_score/.style={precision=2,column name={ARS},
	},
	columns/calinski_harabasz_score/.style={precision=2,column name={CHC},
	},
	columns/davies_bouldin_score/.style={precision=2,column name={DBI},
	},
	columns/silhouette_score/.style={precision=2,column name={$S_{score}$},
	},
    row truelabel/.style={
    every row no #1/.style={
    before row={\rowcolor[gray]{0.9}}},
    },
	% highlight best scores for each measurements
	row AMIS/.style={
	every row #1 column 3/.style={
	% postproc cell content/.append style={
	    % /pgfplots/table/@cell content/.add={\cellcolor{red!10!white}}{}}
     postproc cell content/.style={
          @cell content/.add={$\bf}{$}}
	}
	},
	row ARS/.style={
	every row #1 column 4/.style={
	% postproc cell content/.append style={
	    % /pgfplots/table/@cell content/.add={\cellcolor{red!10!white}}{}}
    postproc cell content/.style={
          @cell content/.add={$\bf}{$}}	
 }
	},
	row h/.style={
	every row #1 column 5/.style={
	% postproc cell content/.append style={
	    % /pgfplots/table/@cell content/.add={\cellcolor{red!10!white}}{}}
	postproc cell content/.style={
          @cell content/.add={$\bf}{$}}
    }
	},
	row c/.style={
	every row #1 column 6/.style={
	% postproc cell content/.append style={
	    % /pgfplots/table/@cell content/.add={\cellcolor{red!10!white}}{}}
    postproc cell content/.style={
          @cell content/.add={$\bf}{$}}	
 }
	},
	row NMIS/.style={
	every row #1 column 7/.style={
	% postproc cell content/.append style={
	    % /pgfplots/table/@cell content/.add={\cellcolor{red!10!white}}{}}
    postproc cell content/.style={
          @cell content/.add={$\bf}{$}}	
 }
	},
	row FMS/.style={
	every row #1 column 8/.style={
	% postproc cell content/.append style={
	    % /pgfplots/table/@cell content/.add={\cellcolor{red!10!white}}{}}
    postproc cell content/.style={
          @cell content/.add={$\bf}{$}}	
    }
	},
	row CHC/.style={
	every row #1 column 9/.style={
	  % postproc cell content/.append style={
            % /pgfplots/table/@cell content/.add={\cellcolor{red!10!white}}{} }
	 postproc cell content/.style={
          @cell content/.add={$\bf}{$}}   
     }
	},
	row DBI/.style={
	every row #1 column 10/.style={
	  % postproc cell content/.append style={
            % /pgfplots/table/@cell content/.add={\cellcolor{red!10!white}}{} }
	 postproc cell content/.style={
          @cell content/.add={$\bf}{$}}   
     }
	},
	row s/.style={
	every row #1 column 11/.style={
	  % postproc cell content/.append style={
            % /pgfplots/table/@cell content/.add={\cellcolor{red!10!white}}{} }
	 postproc cell content/.style={
          @cell content/.add={$\bf}{$}}   
     }
	},
	skip rows between index={0}{9},
	skip rows between index={26}{27},
	row truelabel/.list={0,6,12,17,21,25,29,34,38,43,48,53},
	row AMIS/.list={14,20,25,29,34,38,41,47,52,55,60,68},
	row ARS/.list={14,16,25,29,32,38,41,47,52,55,60,68},
	row h/.list={14,20,25,29,34,38,43,47,52,55,62,68},
	row c/.list={14,16,25,29,32,38,40,47,52,55,59,68},
	row NMIS/.list={14,20,25,29,34,38,41,47,52,55,60,68},
	row FMS/.list={10,17,25,28,32,36,41,45,49,55,59,,68},
	row CHC/.list={10,16,23,30,32,36,41,45,49,57,59,64},
	row DBI/.list={10,20,25,29,32,36,40,46,52,56,59,64},
	row s/.list={10,16,25,29,32,38,40,45,49,55,59,68},
	column addvline/.style={
	every col no #1/.style={
	column type/.add={}{|}}
    },
    column type=r,
	column addvline/.list={2,8},
	every head row/.style={
	before row={
	\caption{Clustering evaluations for real data} \label{tab:real_data_eval_clust} \\
	\toprule
	&  &   & \multicolumn{6}{c|}{extrinsic measurements}& \multicolumn{3}{c}{intrinsic measurements}\\
	},
	after row={
% 	&  &  &  & (rounded) &  (\si{\second}) & (rounded) & (\si{\second}) & \% &  &  &  
	\endfirsthead
	\midrule}
	},
	every last row/.style={after row=\bottomrule}
	]{\data} % filename/path to file
	% }

\pgfplotstableread[col sep=comma]{artificial_data_summary_all_label_clust_eval.tex}\data
	\pgfplotstabletypeset[
	font=\footnotesize,
	multicolumn names, % allows to have multicolumn names
	columns={dataname,k, opt_sos_sdp,adjusted_mutual_info_score,adjusted_rand_score,homogeneity_score,completeness_score,normalized_mutual_info_score,fowlkes_mallows_score,calinski_harabasz_score,davies_bouldin_score,silhouette_score},
	fixed zerofill,
	fixed,
 	columns/dataname/.style={string type,
		column name={instances}
	},
	columns/k/.style={int detect,
		column name={$k$},
	},  
   columns/opt_sos_sdp/.style={precision=0,
        column name = {$d_{SOS}$},
        assign column name/.style={/pgfplots/table/column name={\multicolumn{1}{c|}{##1}}}
    },
	columns/homogeneity_score/.style={precision=2,
	column name={$h$},	},
	columns/completeness_score/.style={precision=2,column name={$c$},
	},
	columns/normalized_mutual_info_score/.style={precision=2,
	column name={NMI},
	},
	columns/fowlkes_mallows_score/.style={precision=2,
	column name={FMS},
	 assign column name/.style={
		/pgfplots/table/column name={\multicolumn{1}{c|}{##1}}}
	},
	columns/adjusted_mutual_info_score/.style={precision=2,column name={AMI},
	},
	columns/adjusted_rand_score/.style={precision=2,column name={ARS},
	},
	columns/calinski_harabasz_score/.style={precision=2,column name={CHC},
	},
	columns/davies_bouldin_score/.style={precision=2,column name={DBI},
	},
	columns/silhouette_score/.style={precision=2,column name={$S_{score}$},
	},
	row AMIS/.style={
	every row #1 column 3/.style={
	% postproc cell content/.append style={
	    % /pgfplots/table/@cell content/.add={\cellcolor{red!10!white}}{}}
    postproc cell content/.style={
          @cell content/.add={$\bf}{$}}	
 }
	},
	row ARS/.style={
	every row #1 column 4/.style={
	% postproc cell content/.append style={
	    % /pgfplots/table/@cell content/.add={\cellcolor{red!10!white}}{}}
    postproc cell content/.style={
          @cell content/.add={$\bf}{$}}	
 }
	},
	row c/.style={
	every row #1 column 5/.style={
	% postproc cell content/.append style={
	    % /pgfplots/table/@cell content/.add={\cellcolor{red!10!white}}{}}
    postproc cell content/.style={
          @cell content/.add={$\bf}{$}}	
 }
	},
	row h/.style={
	every row #1 column 6/.style={
	% postproc cell content/.append style={
	    % /pgfplots/table/@cell content/.add={\cellcolor{red!10!white}}{}}
     postproc cell content/.style={
          @cell content/.add={$\bf}{$}}
	}
	},
	row NMIS/.style={
	every row #1 column 7/.style={
	% postproc cell content/.append style={
	    % /pgfplots/table/@cell content/.add={\cellcolor{red!10!white}}{}}
	postproc cell content/.style={
          @cell content/.add={$\bf}{$}}
    }
	},
	row FMS/.style={
	every row #1 column 8/.style={
	% postproc cell content/.append style={
	    % /pgfplots/table/@cell content/.add={\cellcolor{red!10!white}}{}}
    postproc cell content/.style={
          @cell content/.add={$\bf}{$}}	
 }
	},
	row CHC/.style={
	every row #1 column 9/.style={
	  % postproc cell content/.append style={
            % /pgfplots/table/@cell content/.add={\cellcolor{red!10!white}}{} }
	 postproc cell content/.style={
          @cell content/.add={$\bf}{$}}   
     }
	},
	row DBI/.style={
	every row #1 column 10/.style={
	  % postproc cell content/.append style={
            % /pgfplots/table/@cell content/.add={\cellcolor{red!10!white}}{} }
            postproc cell content/.style={
          @cell content/.add={$\bf}{$}}
	    }
	},
	row s/.style={
	every row #1 column 11/.style={
	  % postproc cell content/.append style={
            % /pgfplots/table/@cell content/.add={\cellcolor{red!10!white}}{} 
    postproc cell content/.style={
          @cell content/.add={$\bf}{$}}
            % }
	    }
	},
    row AMIS/.list={1,7,12,18,20,26,31,36,41,46,51,56},
	row ARS/.list={4,7,12,16,20,25,29,33,39,46,51,56},
	row c/.list={4,7,12,18,20,24,29,36,38,46,48,53},
	row h/.list={4,9,12,18,20,27,31,36,41,46,51,56},
	row NMIS/.list={4,7,12,18,20,25,31,36,39,46,51,56,},
	row FMS/.list={4,7,12,16,20,25,29,33,39,46,51,56},
	row CHC/.list={4,9,12,18,20,25,31,36,38,46,51,55},
	row DBI/.list={3,9,12,18,20,24,30,35,39,46,50,54},
	row s/.list={2,6,12,18,20,24,29,36,39,46,51,55},
	rows truelabel/.style={
    every row no #1/.style={
    before row={\rowcolor[gray]{0.9}}},
    % postproc cell content/.style={
          % @cell content/.add={$\bf}{$}},
    },
    rows truelabel/.list={0,5,10,15,19,23,28,32,37,42,47,52},
	columns addvline/.style={
	every col no #1/.style={
	column type/.add={}{|}}
    },
    column type=r,
	columns addvline/.list={2,8},
	every head row/.style={
	before row={
	\caption{Clustering evaluations for artificial data} \label{tab:artificial_data_eval_clust} \\
	\toprule
	&  & &  \multicolumn{6}{c|}{extrinsic measurements}& \multicolumn{3}{c}{intrinsic measurements}\\
	},
	after row={
% 	&  &  &  & (rounded) &  (\si{\second}) & (rounded) & (\si{\second}) & \% &  &  &  
	\endfirsthead
	\midrule}
	},
	every last row/.style={after row=\bottomrule}
	]{\data} % filename/path to file

We denote by $k_{opt}$ the value of $k$ for which the given quality measure reaches its best value.
We can observe that $k_{opt}$ is not necessarily equal to $k_{true}$.
Table~\ref{table:optimal_K_prop} summarizes the difference $|k_{true}-k_{opt}|$ across all real and artificial datasets, for all quality measures. We can see that ARS has $|k_{true}=k_{opt}|$ for 58.33\% of the artificial datasets and FMS has $|k_{true}=k_{opt}|$ for 50.00\% of the real datasets and for $|k_{true}=k_{opt}|$ for 58.33\% of the artificial datasets. All other measures have $k_{true}\neq k_{opt}|$ for more than 50\% of the datasets.

\begin{table}[htp!]
    \centering
    \footnotesize
    \begin{tabular}{r|rrr|rrr}
    \hline 
     & \multicolumn{3}{c|}{real data} & \multicolumn{3}{c}{artificial data}  \\ 
    & \multicolumn{3}{c|}{$|k_{true}- k_{opt}|$} & \multicolumn{3}{c}{$|k_{true}- k_{opt}|$}  \\ 
    measure & 0 & 1 & 2 & 0 & 1 &2  \\
     & (\%) &  (\%) &  (\%) & (\%) &  (\%) &  (\%)  \\
    \hline 
    AMI        & 25.00 & 8.33  & 66.67 & 25.00 & 16.67 & 58.33 \\
    ARS         & 25.00 & 16.67 & 58.33 & 58.33 & 8.33 &  33.33\\
    $h$         & 16.67 & 25.00 & 58.33 & 16.67 & 0.00 & 83.33 \\
    $c$         & 8.33  & 8.33  & 83.33 & 33.33 & 33.33  & 33.33  \\ 
    NMI         & 8.33  & 8.33  & 66.67 & 41.67 & 0.00    & 58.33\\
    FMS         & 50.00 & 25.00 & 25.00 & 58.33 & 8.33 & 33.33 \\\hline 
    CHC         & 41.67 & 16.67 & 41.67 & 25.00 & 16.67 & 58.33 \\
    DBI         & 25.00 & 41.67 & 33.33 & 33.33 & 41.67 & 25.00 \\  
    $S_{score}$ & 25.00 & 33.33 & 41.67 & 41.67 & 25.00 & 33.33 \\
    \hline
    \end{tabular}
    \caption{Summary of optimal $k_{opt}$ with $|k_{true}-k_{opt}|\in \{0,1,2\}$}
    \label{table:optimal_K_prop}
\end{table}

\subsection{Visualizations}
In this subsection, we visualize the ground truth and optimum clusterings for four datasets: the iris and the sonar datasets from the real dataset group and the 3-spiral and the gaussians1 datasets from the artificial dataset group.
For the iris and gaussian1 datasets, the ground truth and optimum clusterings are very similar, which we can realize by looking at the quality measures for $k=k_{true}$ in the corresponding tables.
For the other two datasets, the quality measures are very poor, which means that the ground truth and the optimum clustering are very different.
For the iris and sonar datasets, which have 4 and 60 dimensions, respectively, only the first two principal components are visualized, while the 3-spiral and gaussians1 datasets are two-dimensional, so the visualizations are completely correct.
We can see that the ground truth clusters for the iris and the Gaussian1 are in the form of ellipsoids that are well separated, so the optimum and the ground clusterings are very similar. For the other two datasets, the ground truth clusters overlap or have the shape of spirals, and optimum clustering is completely different.  

\begin{figure}[htp!]
    \centering
    \subfloat[PCA visualization of ground truth clusters for the \\ iris dataset]{%
        \includegraphics[width=0.4\textwidth]{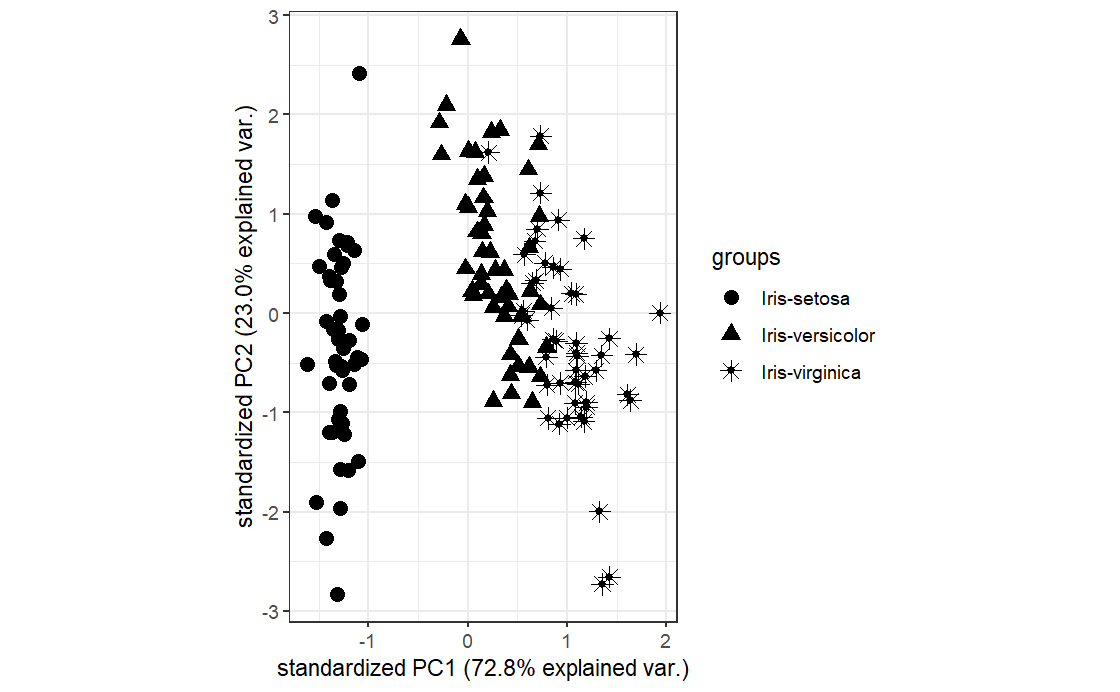}%
        \label{fig:1a}%
        }%
    \hfill%  
    \subfloat[PCA visualization of optimum clusters for $k=k_{true}$ \\ for the iris dataset]{%
        \includegraphics[width=0.4\textwidth]{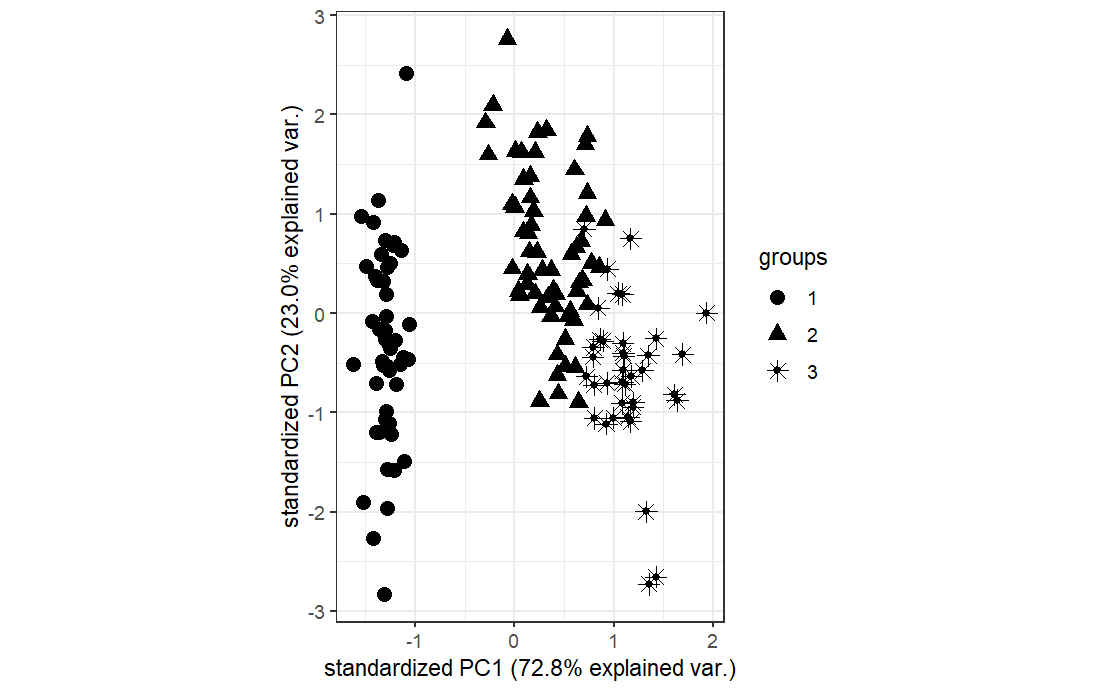}%
        \label{fig:1b}%
        }\\ %
%    \caption{iris dataset}
%\end{figure}
%
%\begin{figure}[H] 
%    \centering
    \subfloat[PCA visualization of ground truth clusters for the \\ sonar  dataset]{%
        \includegraphics[width=0.4\textwidth]{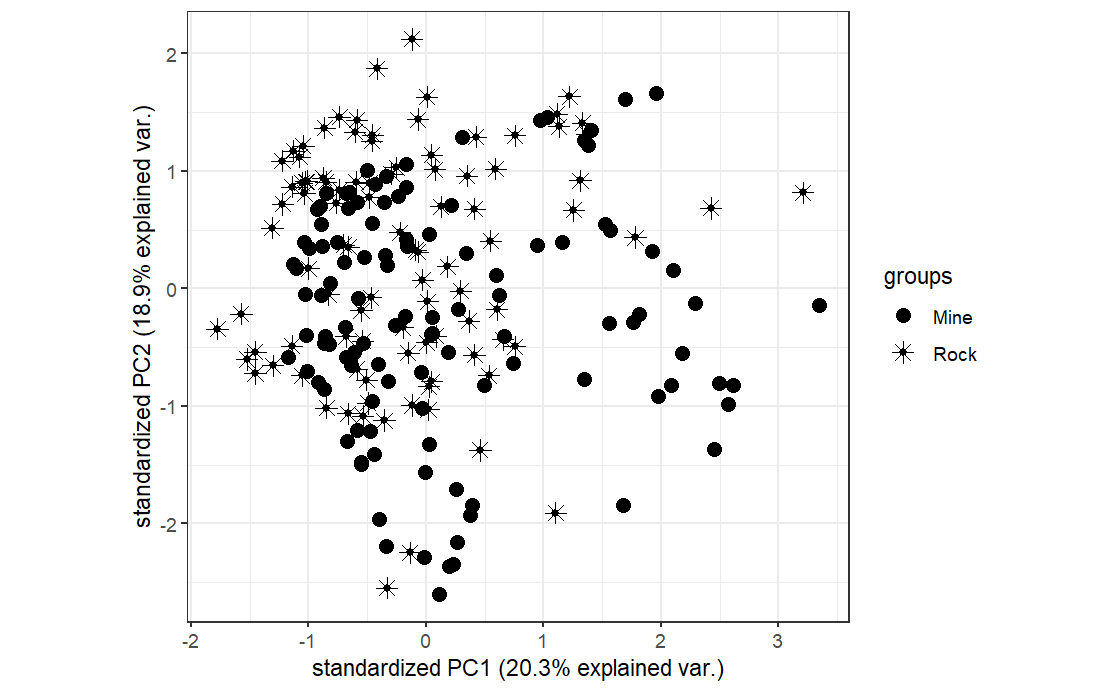}%
        \label{fig:2a}%
        }%
    \hfill%
    \subfloat[PCA visualization of optimum  clusters for $k=k_{true}$ \\ for the sonar dataset]{%
        \includegraphics[width=0.4\textwidth]{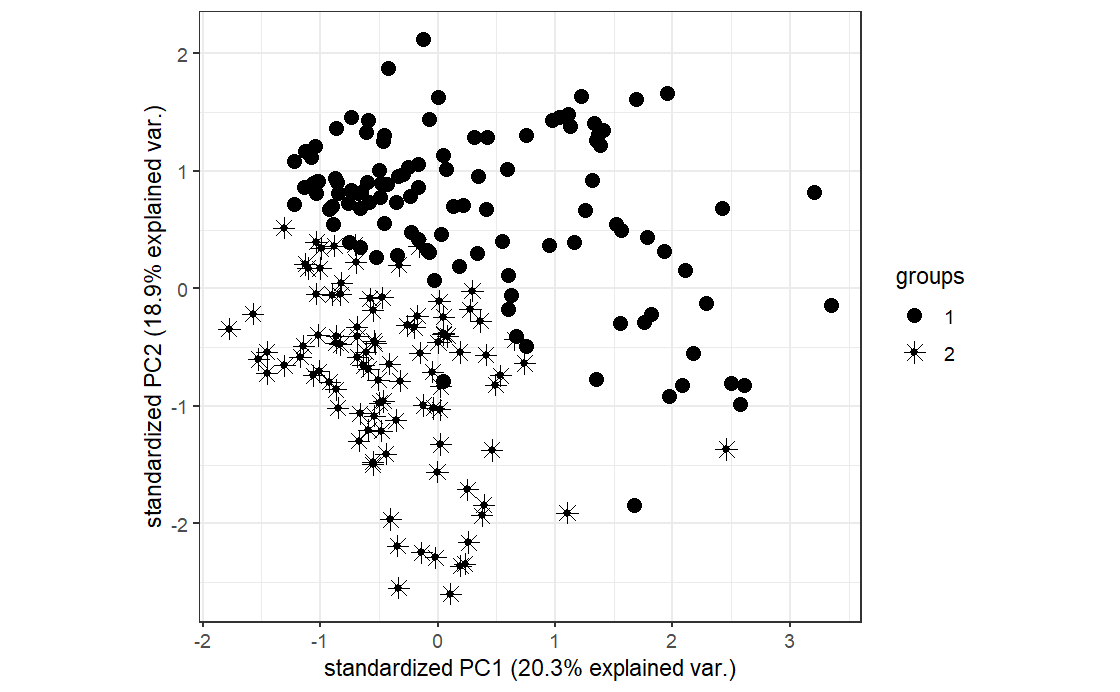}%
        \label{fig:2b}%
        } \\ %
            \caption{Visualizations of ground truth and optimum clusterings for two real datasets with highest and lowest values of quality measures, respectively. }
\end{figure}

\begin{figure}[htp!]
    \centering

    \subfloat[Visualization of ground truth clusters for 3-spiral dataset ]{%
        \includegraphics[width=0.4\textwidth]{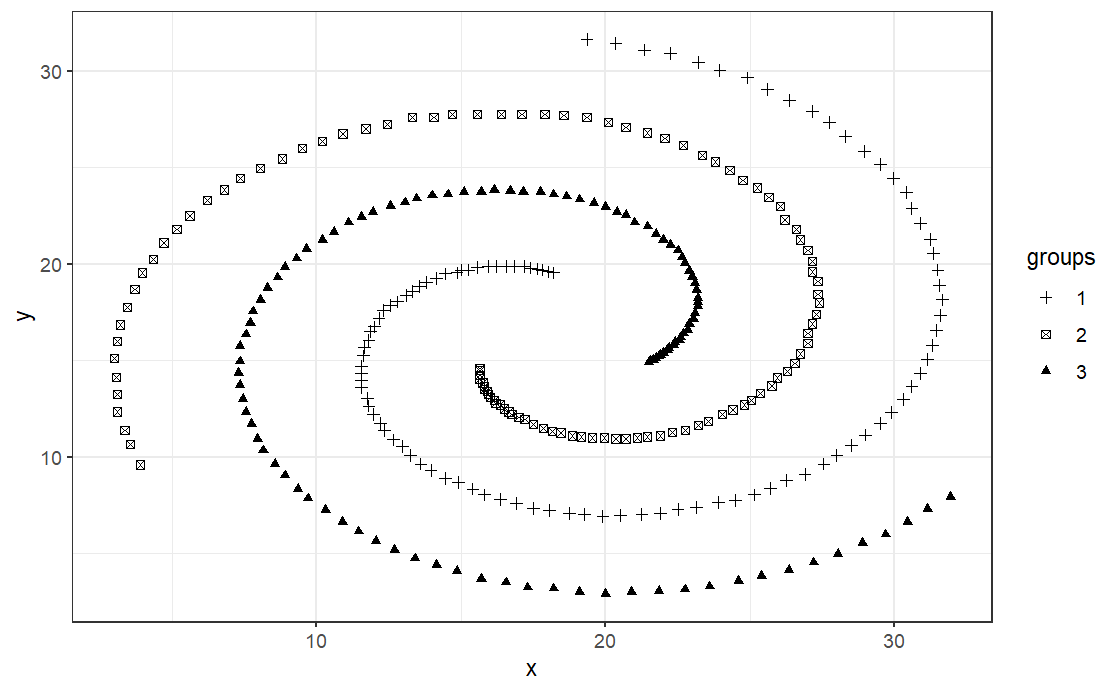}%
        \label{fig:3a}%
        }%
    \hfill%
    \subfloat[Visualization of optimum  clusters for $k=k_{true}$ for \\ 3-spiral dataset]{%
        \includegraphics[width=0.4\textwidth]{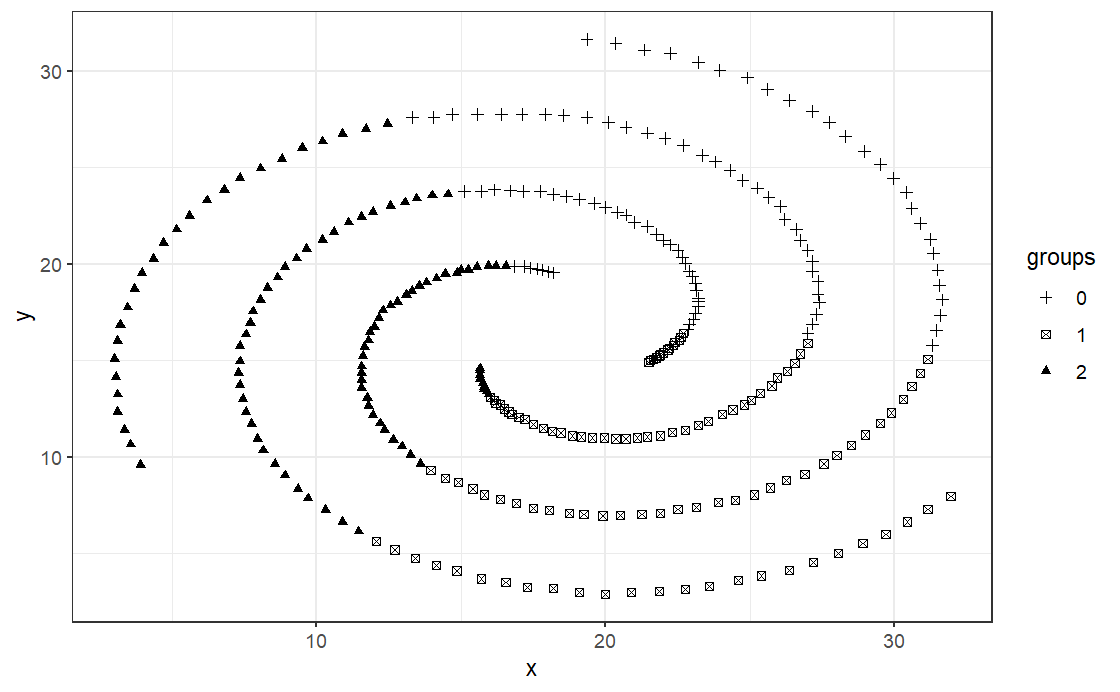}%
        \label{fig:3b}%
        }\\ %

   \centering
    \subfloat[Visualization of ground truth clusters for\\  gaussians1 dataset]{%
        \includegraphics[width=0.4\textwidth]{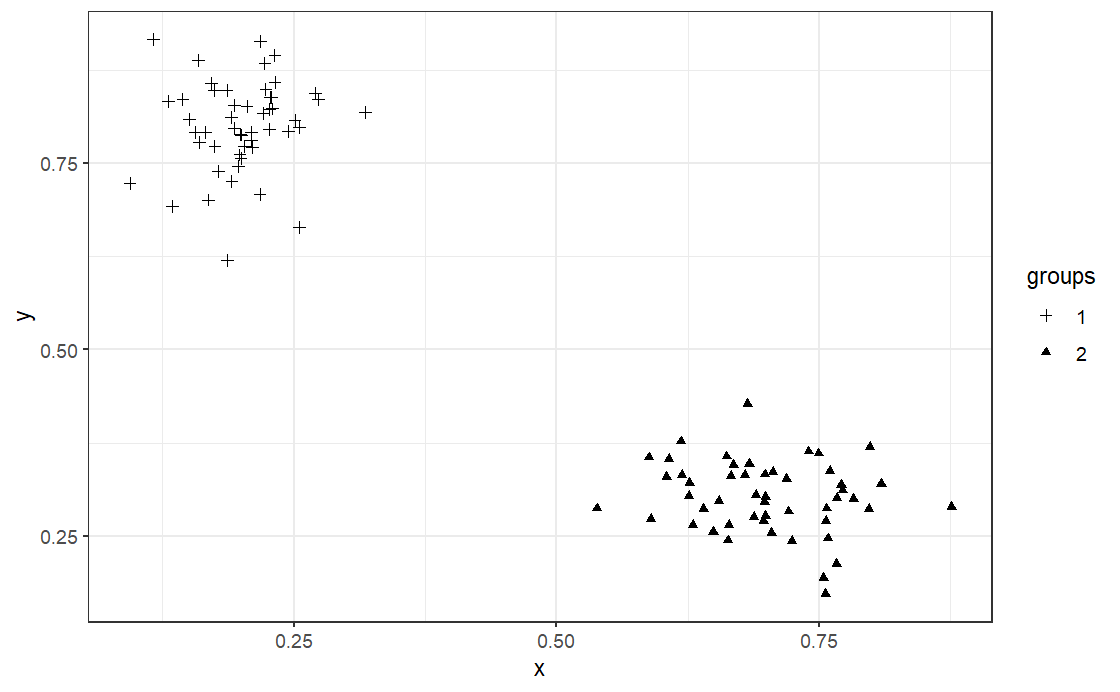}%
        \label{fig:4a}%
        }%
    \hfill%
   \centering
    \subfloat[Visualization of optimum  clusters for $k=k_{true}$ \\ for gaussians1 dataset]{%
        \includegraphics[width=0.4\textwidth]{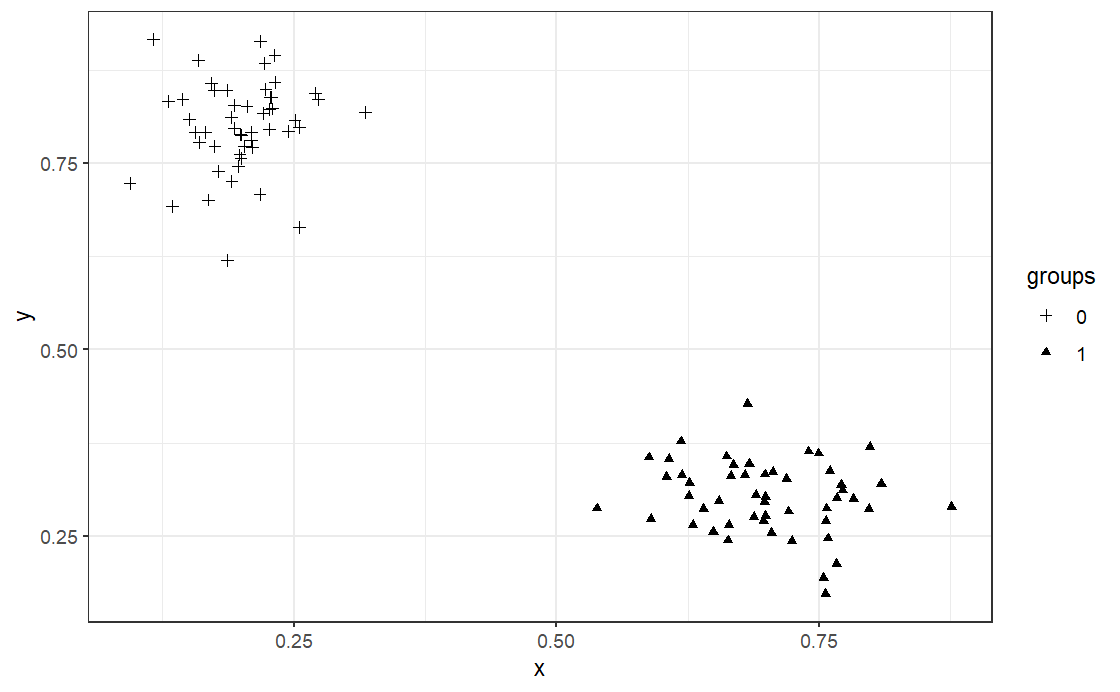}%
        \label{fig:4b}%
        }%
    \caption{Visualizations of ground truth and optimum clusterings for two artificial datasets with highest and lowest values of quality measures, respectively.}
\end{figure}
\clearpage
\section{Discussion}
\label{section:discussion}
The numerical results from Tables \ref{tab:real_data_eval_clust} and \ref{tab:artificial_data_eval_clust} reveal several interesting facts. The main one is that the ground truth clustering usually does not match the optimum MSSC clustering. This has already been observed by other authors \cite{zhang2019critical}, but our work shows it very clearly.

As first,  even detecting $k_{true}$ is very challenging. 
We expected that the most popular quality measures for clustering from Section \ref{sec:quality_measures}  would have the best value for $k=k_{true}$. This  would justify the usual approach where we detect  $k_{true}$ by computing clusterings for different $k$ by using some of the popular heuristic algorithms and then   estimate  $k_{true}$ to be the value for which some of the  quality measures have the best value.

Table \ref{table:optimal_K_prop} shows that this is wrong. None of the  quality measures  can be used to detect $k_{true}$ for MSSC  for the great majority of datasets. The FMS has the best value  for $k=k_{true}$ for 50.0\% of real datasets and for 58.3\% of artificial datasets and ARS has the best value for 58.3\% of artificial datasets.
All the other measures achieve  best values for $k\neq k_{true}$ in more than 50 \%. 
This  means that we can not rely on these quality measures  to detect $k_{true}$.
Also the optimum value of \eqref{eq:mssc} can not be used to detect $k_{true}$ since this value is monotonically decreasing with $k$ increasing, as is depicted in the  3rd columns of Tables \ref{tab:real_data_eval_clust} and \ref{tab:artificial_data_eval_clust}. 
It is also interesting that the objective value of \eqref{eq:mssc}, which is obtained on the ground truth clustering, is usually far from the optimum value for $k=k_{true}$, which means that the ground truth clustering is a feasible solution for \eqref{eq:mssc}, which is in most cases far from the optimum.

We note that the intrinsic measures need special attention: all of them are related to Euclidean distance, but none of them measures exactly the same thing as we do in the objective function of \eqref{eq:mssc}. They compare the within-the-group  distances to between-the-group  distances in  different ways, whereas the objective function of \eqref{eq:mssc} only considers within-the-group distances. This is also a possible reason why these measures usually have the best value at different $k$, also different from $k_{true}$.

The only exceptions are the gaussians1 dataset, where the ground truth and the optimum clustering are the same, and the iris dataset, where the ground truth and the optimum clustering are similar but not the same.
These clusterings are visualized in Figures \ref{fig:1a}--\ref{fig:1b} and \ref{fig:4a}--\ref{fig:4b}.
These figures actually describe the numerical results very well. If the ground truth clusters have the expected  geometry, i.e.,  the clusters have the form of convex sets, ellipsoids, which are well separated from each other, then such ground truth clustering is  very similar to optimum clustering.
Otherwise, if the clusters have geometrically different shapes, such as 3-spirals, or are overlapping, then the ground truth clustering is far from the optimum clustering, which always enforces the ellipsoidal  geometry of the clusters since the Euclidean distance underlies this model.

Our work therefore confirms the conclusions of \cite{zhang2019critical} that the datasets with known ground truth, which are usually used as benchmark datasets for classification problems, can be used for benchmarking the clustering algorithms with careful attention, since the class labels are usually assigned based on the properties of each individual data point, probably including additional information not present in the dataset itself, while the clustering algorithms take into account the relationships between the data.

\section{Conclusions}

In this paper, we considered the mathematical programming formulation of the NP -hard minimum sum-of-squares clustering problem (MSSC) and solved it to optimality for a number of real and artificial datasets for which the ground truth clustering (the clustering created by the data provider) was available and which were of small size (the number of data points times the number of true clusters $k_{true}$ had to be approximately less than 1000). We solved these instances for the number of clusters $k$ equal or close to the ground truth $k_{true}$ using the exact solver SOS-SDP \cite{piccialli2022sos}.

 For each dataset and each $k$ that we used we compared the optimum clustering (optimum solution of \eqref{eq:mssc})
with  the ground truth clustering by using a number of extrinsic and intrinsic measures.

We  showed that the ground truth clusterings are usually quite far from the optimum clustering, for all $k$ close to $k_{true}$, which means that (i) they yield the value of the objective function of \eqref{eq:mssc}, which is usually much worse (higher) compared to the optimum value of \eqref{eq:mssc}, (ii) the values of intrinsic measures evaluated at the ground truth clustering are usually much worse than the values of the intrinsic measures evaluated at the optimum clustering, (iii) the values of the extrinsic measures which we used to measure the alignment between the ground truth and the optimum clustering showed that they differ a lot.
However, if the ground truth clustering has a natural geometry, i.e., if the clusters look like ellipsoids that are well separated from each other, then the ground truth clustering and the optimum clustering are very similar.

We can derive the following main conclusions: (i) The ground truth clusterings were defined by data providers who often used similarity measures that were not based on Euclidean distance or on distances equivalent to this distance. It is  likely that the similarity measure used was in fact not a distance (metric) according to the mathematical definition, see e.g. \cite{sharma2006metric}.
They may have even used additional information not included in the variables describing the data points.
 Therefore, there is most likely no mathematical distance at which the optimal value of \eqref{eq:mssc}, where the objective function would be defined using this distance, would have ground truth clustering as the optimum solution.
 (ii) We should be very careful when comparing the clustering obtained by a particular clustering algorithm with the ground truth clustering. If such an algorithm measures the similarity between data points with a distance equal to the Euclidean distance, such as the famous $k$-means algorithm, while the ground truth clusters do not have an ellipsoidal geometry, then we cannot expect to obtain a solution that is close to the ground truth clustering.
 (iii) Determining the most appropriate number of clusters by considering where the values of the extrinsic or intrinsic measures have the best value can also be misleading. Very often these measures give contradictory  answers: different measures suggest different $k$, which very often differ from $k_{true}$. The situation is somewhat better when the clusters have a natural ellipsoidal geometry.

This confirms the great importance of choosing the similarity measure over the data points in data clustering: it should capture the geometry of the underlying data. Clustering should also be evaluated using quality measures that are aligned with the similarity measure used in the calculation.

\bibliographystyle{apalike}       % APS-like style for physics
\bibliography{References}   % name your BibTeX data base

\end{document}